\documentclass[12pt]{article}
\usepackage{amsfonts}
\usepackage{latexsym, amssymb, amsmath, amscd, amsfonts, epsfig, graphicx, colordvi}
\usepackage{ifpdf}
\usepackage{graphicx}
\usepackage{pstricks}
\newtheorem{thm}{Theorem}[section]

\newtheorem{defi}[thm]{Definition}
\newtheorem{lem}[thm]{Lemma}

\def\pf{\noindent{\it Proof.} }
\def\qed{\nopagebreak\hfill{\rule{4pt}{7pt}}
\medbreak}

\def\qed{\nopagebreak\hfill{\rule{4pt}{7pt}}
\medbreak}

\begin{document}

\begin{center}
{\bf \large A Unification of Two Refinements of Euler's Partition
Theorem}
\end{center}

\begin{center}
 William Y. C. Chen$^1$, Henry Y. Gao$^2$,  Kathy Q. Ji$^3$ and Martin Y. X.
 Li$^4$

   Center for Combinatorics, LPMC-TJKLC\\
   Nankai University\\
    Tianjin 300071, P.R. China

   \vskip 1mm

Email: $^1$chen@nankai.edu.cn, $^2$gaoyong@cfc.nankai.edu.cn,
$^3$ji@nankai.edu.cn, $^4$lyz6988@yahoo.com.cn

\vskip 3mm

Dedicated to Professor George Andrews on the Occasion of His Seventieth
 Birthday
\end{center}

\vskip 6mm \noindent{{\bf Abstract.} We obtain  a unification of two
refinements of Euler's partition theorem respectively due to Bessenrodt and
Glaisher. A specialization of Bessenrodt's insertion algorithm for a
generalization of the Andrews-Olsson partition identity is used
in our combinatorial construction.



\section{ Introduction}

There are several  bijective proofs and refinements of the classical
partition theorem of Euler. This paper will be concerned with two
remarkable bijections obtained by Sylvester \cite{syl82} and
Glaisher \cite{glaisher}, see also, \cite[pp.8--9]{and04}.
Glaisher's bijection implies  a refinement of Euler's theorem
involving the number of odd parts in a partition with distinct parts
and the number of parts repeated odd times in a partition with odd
parts. On the other hand, as observed by Bessenrodt \cite{bes94},
Sylvester's bijection also leads to a refinement of Euler's theorem.
The main result of this paper is a unification of these two
refinements that does not directly follow from Sylvester's bijection
and Glaisher's bijection.

Let us give an overview of the background and  terminology. We will
adopt the common notation on partitions used in Andrews
 \cite[Chapter 1]{and76}.
{\it A partition} $\lambda$ of a positive integer $n$ is a finite
nonincreasing sequence of positive integers
$$\lambda=(\lambda_1,\,\lambda_2,\ldots,\,\lambda_r)$$ such that
$\sum_{i=1}^r\lambda_i=n.$  The entries $\lambda_i$ are called the
parts of $\lambda$, and $\lambda_1$ is the largest part. The
number of parts of  $\lambda$ is called the length of $\lambda$,
denoted by $l(\lambda)$. The weight of $\lambda$ is the sum of its
parts, denoted by $|\lambda|$.  A partition $\lambda$ can also be
represented in the following form \[
\lambda=(1^{m_1},\,2^{m_2},\,3^{m_3}, \ldots),\] where
  $m_i$ is the multiplicity of the part $i$ in $\lambda$. The
conjugate partition of $\lambda$ is defined by
$\lambda'=(\lambda_1',\lambda_2',\ldots,\,\lambda'_t)$, where
$\lambda'_i$ is the number of parts of $\lambda$ that are greater
than or equal to $i$.

Euler's partition theorem reads as follows.

\begin{thm}[Euler]\label{euler}
The number of partitions of $n$ with distinct parts is equal to the
number of partitions of $n$ with odd parts.
\end{thm}

Let $\mathcal{D}$  denote the set of partitions with distinct
parts, and let $\mathcal {D}(n)$ denote the set of partitions of
$n$ in $\mathcal {D}$. Similarly, let $\mathcal{O}$  denote the
set of partitions with odd parts, and  let $\mathcal {O}(n)$
denote the set of partitions of $n$ in $\mathcal {O}$.
 Sylvester's fish-hook
bijection  \cite{syl82}, also referred to as Sylvester's
bijection, and Glaisher's bijection \cite[pp.8-9]{and04} have
established direct correspondences between $\mathcal{D}(n)$  and
$\mathcal{O}(n)$. These two bijections imply refinements of
Euler's theorem. There are also several other refinements of
Euler's partition theorem, see, for example,
\cite{ald69,and66-a,And83,bes94,kim99,pak03,zen05},
\cite[pp.51--52]{bre99}, \cite[pp.46--47]{fin88}.

Sylvester's refinement \cite[p.24]{and76} is stated as follows.
Recall that a chain in a partition with distinct parts
 is a maximal sequence
of parts consisting of consecutive integers. The number of chains in
a partition $\lambda$ is denoted by $n_c(\lambda)$. The number of
different parts in a partition $\mu$ is denoted by $n_d(\mu)$.  For
example, the partition $(8,7, 5, 3,2,1)$ has three chains, and the
partition $(8,6,6,5,4,4,2,1)$ has  six different parts.

\begin{thm}[Sylvester]\label{syl}
The number of partitions of $n$ into distinct parts with exactly $k$
chains   is equal to  the number of partitions of $n$ into odd parts
{\rm(}repetitions allowed{\rm)} with exactly $k$ different parts. In
the notation of generating functions, we have
\begin{equation}
\sum_{\lambda\in \mathcal{D}}z^{n_c(\lambda)}q^{|\lambda|}=\sum_{\mu
\in \mathcal{O}}z^{n_d(\mu)}q^{|\mu|}.
\end{equation}
\end{thm}

Fine \cite[pp.46--47]{fin88} has derived  a refinement of Euler's
theorem.

\begin{thm}[Fine]\label{fine1} The  number of partitions of $n$ into distinct parts with
largest part $k$ is equal to the number of partitions of $n$ into
odd parts such that  the largest part plus twice the number of parts
equals $2k+1$. In the notation of generating functions, we have
\begin{equation} \label{fine1-e}
\sum_{\lambda\in \mathcal{D}}x^{\lambda_1}q^{|\lambda|}=\sum_{\mu
\in \mathcal{O}}x^{(\mu_1-1)/2+l(\mu)}q^{|\mu|}.
\end{equation}
\end{thm}

 Bessenrodt
\cite{bes94} has shown that Sylvester's bijection implies the
following refinement,
 which is a limiting
 case of the lecture hall theorem due to  Bousquet-M\'{e}lou
and Erikssonin \cite{bou97, bou972}. Let $l_a(\lambda)$ denote the
alternating sum of $\lambda$, namely,
\[ l_a(\lambda)=\lambda_1-\lambda_2+\lambda_3-\lambda_4+\cdots.\]

\begin{thm}[Bessenrodt]\label{euler-bousquet}The
number of partitions of $n$ into distinct parts with alternating sum
$l$ is equal to the number of partitions of $n$ with $l$ odd parts.
In terms of generating functions, we have
\begin{equation}\label{euler-bousquet-e}
\sum_{\lambda \in \mathcal{D}}
y^{l_a(\lambda)}q^{|\lambda|}=\sum_{\mu \in
\mathcal{O}}y^{l(\mu)}q^{|\mu|}.
\end{equation}
\end{thm}

It has also been shown by Bessenrodt \cite{bes94} that Sylvester's
bijection maps the parameter $n_c(\lambda)$ to the parameter
$n_d(\mu)$. Combining the above Theorems \ref{syl} and \ref{fine1},
we arrive at the following equidistribution result.

\begin{thm}[Sylvester-Bessenrodt]
The  number of partitions of $n$ into distinct parts with largest
part $k$, alternating sum $l$ and $m$ chains
 is equal to the number of partitions of $n$ into $l$ odd parts
with exactly $m$ different parts such that the largest part plus
twice the number of parts equals $2k+1$. In terms of generating
functions, we have
\begin{equation}\label{Sylvester-Bessenrodt-e}
\sum_{\lambda \in \mathcal{D}}x^{\lambda_1}
y^{l_a(\lambda)}z^{n_c(\lambda)}q^{|\lambda|} =\sum_{\mu \in
\mathcal{O}}x^{(\mu_1-1)/2+l(\mu)}y^{l(\mu)}z^{n_d(\mu)}q^{|\mu|}.
\end{equation}
\end{thm}

Recently, Zeng \cite{zen05} has found a generating function proof of
the above three-parameter  refinement
\eqref{Sylvester-Bessenrodt-e}.

From a different angle,
 Glaisher \cite{glaisher}, see also \cite[pp.8--9]{and04}, has given
  a refinement of
 Euler's partition theorem. Let
 $l_o(\lambda)$ denote the number of odd parts in $\lambda$,
and let $n_o(\mu)$ denote the number of different parts in $\mu$
with odd multiplicities.

\begin{thm}[Glaisher]\label{euler-Glaisher}
The number of partitions of $n$ into distinct parts with $k$ odd
parts is equal to the number of partitions of $n$ with odd parts
such that there are exactly $k$ different parts repeated odd times.
In terms of generating functions, we have
\begin{equation}\label{euler-Glaisher-e}
\sum_{\lambda \in
\mathcal{D}}x^{l_o(\lambda)}q^{|\lambda|}=\sum_{\mu \in
\mathcal{O}}x^{n_o(\mu)}q^{|\mu|}.
\end{equation}
\end{thm}

Given the two bijections of Sylvester and Glaisher, it is natural to
ask the question whether the joint distribution of the statistics
$(l_o(\lambda), l_a(\lambda))$ of partitions of
 $n$ with distinct
parts coincides with the joint distribution of the statistics
$(n_o(\mu), l(\mu))$ of partitions with odd parts. It turns out that
this is indeed the case. However, neither  Sylvester's bijection nor
Glaisher's bijection implies this result. To give a combinatorial
proof of this result, we need Bessenrodt's insertion algorithm.

It should be noted that the equidistriubtion of $(l_o(\lambda),
l_a(\lambda))$  and $(n_o(\mu), l(\mu))$ can also be deduced from a
recent result of Boulet \cite{bou} by the manipulation of generating
functions.

This paper is organized as follows. In Section 2, we present the
main result and some  lemmas.  Section 3 is devoted to a brief
review of Bessenrodt's insertion algorithm. In Section 4, we utilize
Boulet's formula to give
 a generating function proof of the
two-parameter refinement of Euler's theorem. In Section 5, we give
a combinatorial proof of the unification of the refinements of
Bessenrodt \eqref{euler-bousquet-e} and Glaisher
\eqref{euler-Glaisher-e}.

\setcounter{equation}{0}
\section{The main result}

The main result of this paper is
the following unification of the refinements  of
Bessenrodt and Glaisher.

\begin{thm}\label{main2}
The number of partitions of $n$ into distinct parts with $l$ odd
parts and alternating sum $m$ is equal to the number of partitions
of $n$ into exactly $m$ odd parts and $l$ parts repeated odd times.
In terms of generating functions, we have
\begin{equation}\label{main2-e}
\sum_{\lambda \in
\mathcal{D}}x^{l_o(\lambda)}y^{l_a(\lambda)}q^{|\lambda|} =\sum_{\mu
\in \mathcal{O}}x^{n_o(\mu)}y^{l(\mu)}q^{|\mu|}.
\end{equation}

\end{thm}

For example, Table \ref{tab1} illustrates  the case of $n=7$.

\begin{table}[htb]
    \begin{center}
    \begin{tabular}{|c|c|c||c|c|c|c|c}
         \hline
    {\small $\lambda \in \mathcal{D}(7)$} &
    {\small  $l_o(\lambda)$} &
    {\small $l_a(\lambda)$}& {\small $\mu \in \mathcal{O}(7)$} &

   {\small  $n_o(\mu)$} &
    {\small  $l(\mu)$} \\
         \hline
        $(7)$& $1$&$7$ &$(1^7)$&$1$&$7$\\[5pt]

        $(1,6)$&$1$&$5$& $(1^4,3)$&$1$&$5$\\[5pt]

         $(2,5)$&$1$& $3$&$(1,3^2)$&$1$& $3$ \\[5pt]

       $(3,4)$&$1$&$1$&  $(7)$&$1$&$1$\\[5pt]
         $(1,2,4)$&$1$&$3$& $(1^2,5)$&$1$&$3$\\[5pt]\hline
    \end{tabular}
    \end{center}
\caption{The case of $n=7$ for Theorem \ref{main2}.}\label{tab1}
\end{table}

It is clear that the above theorem reduces to
 Bessenrodt's refinement (1.3) when $x=1$ and to Glaisher's
refinement (1.5)  when $y=1$.

To prove Theorem \ref{main2}, we proceed to construct a bijection
$\Delta$ between $\mathcal{D}(n)$ and $\mathcal{O}(n)$ such that for
$\lambda \in \mathcal{D}(n)$ and
$\mu=\Delta(\lambda)\in\mathcal{O}(n)$, we have
\[l_o(\lambda)=n_o(\mu),\quad l_a(\lambda)=l(\mu).\]

Let $\mathcal{A}_1(n)$ denote the set of partitions of $n$ subject
to the following conditions:
\begin{enumerate}
\item Only parts divisible by $2$ may be
repeated.
\item
The difference between successive parts is at most $4$ and strictly
less than $4$ if either part is divisible by $2$.
\item
The smallest part is less than $4$.
\end{enumerate}

 By considering the conjugate of the 2-modular representation
 of a partition, it is easy to establish a bijection between
$\mathcal{D}(n)$ and $\mathcal{A}_1(n)$.

\begin{lem}\label{lem1}
There is a bijection $\varphi$ between $\mathcal{D}(n)$ and
$\mathcal{A}_1(n)$. Furthermore, for $\lambda \in \mathcal{D}(n)$
and $\alpha=\varphi(\lambda)\in \mathcal{A}_1(n)$, we have
\begin{equation}
l_o(\lambda)=l_o(\alpha),\quad
l_a(\lambda)=2r_2(\alpha)+l_o(\alpha),
\end{equation}
where $r_2(\alpha)$ denotes the number of parts congruent to $2$
modulo $4$ in $\alpha$.
\end{lem}

Let $\mathcal{A}_2(n)$ denote the set of partitions of $n$ subject
to the following conditions:
\begin{enumerate}
\item No part divisible by $4$.
\item Only parts divisible by $2$ may be repeated.
\end{enumerate}

We then establish a bijection between $\mathcal{O}(n)$ and
$\mathcal{A}_2(n)$ in the spirit of Glaisher's bijection.

\begin{lem}\label{lem2}
There is a bijection $\psi$ between $\mathcal{O}(n)$ and
$\mathcal{A}_2(n)$. Furthermore, for $\mu \in \mathcal{O}(n)$ and
$\beta=\psi(\mu)\in \mathcal{A}_2(n)$, we have
\begin{equation}
n_o(\mu)=l_o(\beta),\quad  l(\mu)=2r_2(\beta)+l_o(\beta).
\end{equation}
\end{lem}

In view of  the above two lemmas, we see that
 Theorem \ref{main2} can be deduced from
the following theorem.

\begin{thm}\label{temp-t}
There is a bijection $\phi$ between $\mathcal{A}_1(n)$ and
$\mathcal{A}_2(n)$. Furthermore, for $\alpha \in \mathcal{A}_1(n)$
and $\beta=\phi(\alpha)\in \mathcal{A}_2(n)$, we have
\begin{equation}
l_o(\alpha)=l_o(\beta),\quad r_2(\alpha)=r_2(\beta).
\end{equation}
\end{thm}

We find that Theorem \ref{temp-t} can be deduced from Bessenrodt's
insertion algorithm which was devised as a  combinatorial proof of a
generalization of Andrews-Olsson's theorem \cite{bes95}.  Combining
the bijection $\varphi$ for Lemma \ref{lem1}, $\psi$ for Lemma
\ref{lem2} and $\phi$ for Theorem \ref{temp-t}, we are led to a new
bijection $\Delta$ for Euler's partition theorem which implies the
equidistribution of the statistics $(l_o(\lambda), l_a(\lambda))$ of
partitions with distinct parts and the statistics $(n_o(\mu),
l(\mu))$ of partitions with odd parts.

\section{Bessenrodt's insertion algorithm}

To provide a purely combinatorial proof of  Andrews-Olsson's theorem
\cite{And91}, Bessenrodt \cite{bes91} constructs an explicit
bijection on the sets of partitions in Andrews-Olsson's theorem,
which we call Bessenrodt's insertion algorithm. The original
insertion algorithm does not imply the bijection in Theorem
\ref{temp-t}, but we find that the generalized insertion algorithm
given by Bessenrodt \cite{bes95} in 1995 can be used to establish
the bijection required by Theorem \ref{temp-t}.

We give an overview  of Bessenrodt's insertion algorithm.
 Let $N$ be an integer, and let
$\mathbb{A}_N=\{a_1,a_2,\ldots,a_r\}$ with $1\le a_1<
a_2<\cdots<a_r<N$. Andrews-Olsson's theorem involves two sets
$\mathcal {AO}_1(\mathbb{A}_N;n,N)$ and $\mathcal
{AO}_2(\mathbb{A}_N;n,N)$ defined below.

\begin{defi}
Let $\mathcal {AO}_1(\mathbb{A}_N;n,N)$ denote the set of
partitions of $n$ satisfying the following conditions:
\begin{enumerate}
\item Each part is congruent to 0 or some $a_i$ modulo $N$;

\item Only the multiples of $N$ can be repeated;

\item The difference between two successive parts is at most $N$
and strictly less than $N$ if either part is divisible by $N$;

\item The smallest part is less than $N$.
\end{enumerate}
\end{defi}

\begin{defi}
Let $\mathcal {AO}_2(\mathbb{A}_N;n,N)$ denote the set of
partitions of $n$ satisfying the following conditions:
\begin{enumerate}
\item Each part is congruent to some $a_i$ modulo $N;$

\item No part can be repeated.
\end{enumerate}
\end{defi}

The cardinalities of $\mathcal {AO}_1(\mathbb{A}_N;n,N)$ and
$\mathcal {AO}_2(\mathbb{A}_N;n,N)$ are denoted by
 $p_1(\mathbb{A}_N;n,N)$ and $p_2(\mathbb{A}_N;n,N)$ respectively.
Andrews-Olsson's theorem  is  stated as follows.

\begin{thm}[Andrews-Olsson]\label{and-ols}
For any $n \in\mathbb{N}$, we have
$$p_1(\mathbb{A}_N;n,N)=p_2(\mathbb{A}_N;n,N).$$
\end{thm}

By examining the two sets $\mathcal {A}_1(n)$ and $\mathcal
{A}_2(n)$ in Theorem \ref{temp-t}, we find they are somehow
analogous to the  two sets $\mathcal {AO}_1(\mathbb{A}_N;n,N)$ and
$\mathcal {AO}_2(\mathbb{A}_N;n,N)$ in Andrews-Olsson's theorem.
Moreover, we could also apply Bessenrodt's insertion algorithm to
establish a bijection between $\mathcal {A}_1(n)$ and $\mathcal
{A}_2(n)$. Here we present a more general bijection $\Phi$ between
the two sets $\mathcal {C}_1(\mathbb{A}_{2N};n, 2N)$ and $\mathcal
{C}_2(\mathbb{A}_{2N}; n, 2N)$, and we can restrict the bijection
$\Phi$ to $\mathcal {A}_1(n)$ and $\mathcal {A}_2(n)$
 by  setting $N=2$ and $\mathbb{A}_4=\{1,2,3\}$.

\begin{defi}
Let $\mathcal {C}_1(\mathbb{A}_{2N};n, 2N)$ denote the set of
partitions of $n$ satisfying the following conditions:

\begin{enumerate}

\item Each part is congruent to 0 or some $a_i$ modulo $2N$;

\item Only the multiples of  $N$ can be repeated;

\item The difference between two successive parts is at most $2N$
and strictly less than $2N$ if either part is divisible by $N;$

\item The smallest part is less than $2N.$
\end{enumerate}
\end{defi}

\begin{defi}
Let $\mathcal {C}_2(\mathbb{A}_{2N}; n,2N)$ denote the set of
partitions of $n$ satisfying the following conditions:

\begin{enumerate}

\item Each part is congruent to some $a_i$ modulo $2N$;

\item Only multiples of $N$ may be repeated;

\end{enumerate}

\end{defi}

The cardinalities of $\mathcal {C}_1(\mathbb{A}_{2N};n, 2N)$ and
$\mathcal {C}_2(\mathbb{A}_{2N};n, 2N)$ are denoted by
 $c_1(\mathbb{A}_{2N};n, 2N)$ and $c_2(\mathbb{A}_{2N};n, 2N)$ respectively.
Then we have the following theorem which will be needed to prove
Theorem \ref{temp-t}.

\begin{thm}\label{cgeneral}
For any $n \in\mathbb{N}$, we have
$$c_1(\mathbb{A}_{2N};n, 2N)=c_2(\mathbb{A}_{2N};n, 2N).$$
\end{thm}

Theorem \ref{cgeneral} can be proved either by a variant of
Bessenrodt's insertion algorithm obtained in 1991, or by
specializing a  generalization of Bessenrodt's algorithm
obtained in 1995.

We outline the first approach by constructing a
bijection $\Phi$ between $\mathcal {C}_1(\mathbb{A}_{2N};n, 2N)$
and $\mathcal {C}_2(\mathbb{A}_{2N};n, 2N)$ based on a
variant of Bessenrodt's
insertion algorithm \cite{bes91}.

 For $\lambda\in\mathcal
{C}_1(\mathbb{A}_{2N};n, 2N),$ we first extract some parts from
$\lambda$ to form a pair of partitions $(\alpha,\beta),$ where
$\alpha\in\mathcal {C}_1(\mathbb{A}_{2N};n, 2N)\cap\mathcal
{C}_2(\mathbb{A}_{2N};n, 2N)$ and $\beta$ is a partition with
parts divisible by $2N.$ Then we insert $\beta$ into $\alpha$ to
get a partition $\gamma\in\mathcal {C}_2(\mathbb{A}_{2N};n, 2N).$

{\noindent The bijection $\Phi$ consists of the following two
steps.}

{\noindent  Step 1:} Extract certain parts from
$\lambda=(\lambda_1,\lambda_2,\ldots,\lambda_{l(\lambda)})\in\mathcal
{C}_1(\mathbb{A}_{2N};n, 2N).$

 We now construct a pair of partitions
$(\alpha,\beta)$  based on the partition $\lambda$. Let $\lambda
_{j}$ be a part divisible by $2N$, and $\lambda _{t}$ be the
smallest part  bigger than $\lambda _j$. We remove $\lambda _{j}$ if
$\lambda _{t}$ does not exist or the difference between $\lambda
_{t}$ and $\lambda _{j+1}$ satisfies the difference condition to
$\mathcal {C}_1(\mathbb{A}_{2N};n, 2N)$. After removing these parts
$\lambda_j$,  we obtain a partition $\alpha^1$ in $\mathcal
{C}_1(\mathbb{A}_{2N};n, 2N)$, and we can rearrange these parts that
have been removed  to form partition $\beta^1$.

Assume that there are $l$  parts divisible by $2N$ in $\alpha^1$.
Let $t=1$. We may iterate the following procedure  until we get a
pair of partitions $(\alpha^{l+1},\,\beta^{l+1})$.
\begin{itemize}
\item Let $\alpha^t_i$ be the largest part divisible by $2N$ in
$\alpha^t$.

\item Subtract $2N$ from
$\alpha^t_1,\,\alpha^t_2,\ldots,\,\alpha^t_{i-1}$ and  remove
$\alpha^t_i$ from $\alpha^t$.

\item Rearrange the remaining parts to give a new partition
$\alpha^{t+1}$ and add one part  $(i-1)\cdot 2N+\alpha^t_i$ to
$\beta^{t}$ to get $\beta^{t+1}$.
\end{itemize}
Then let $\alpha=\alpha^{l+1},\,\beta=\beta^{l+1} .$  It can be
seen that  $\alpha \in \mathcal {C}_1(\mathbb{A}_{2N};n, 2N)\cap
\mathcal {C}_2(\mathbb{A}_{2N};n, 2N)$ and $\beta_1\leq 2N\cdot
l(\alpha)$.

{\noindent Step 2:} Insert $\beta$ into $\alpha$ to generate a
partition $\gamma\in\mathcal {C}_2(\mathbb{A}_{2N};n, 2N).$

{\noindent}For each $\beta_i,$ we add $2N$ to the first $\beta_i/2N$
parts of $\alpha\colon
\alpha_1,\alpha_2,\ldots,\alpha_{\beta_i/2N},$ then denote the
resulted partition by $\gamma$. It can be shown that
  $\gamma \in \mathcal {C}_2(\mathbb{A}_{2N};n, 2N)$ for $\beta_1\leq 2N\cdot
l(\alpha)$. For the details of the proof, see \cite{bes91,yee02}.

The inverse map $\Phi^{-1}$ can be described as follows.  For
$\gamma\in\mathcal {C}_2(\mathbb{A}_{2N};n, 2N)$, we  extract
certain parts from $\gamma$ to get a pair of partitions
$(\alpha,\beta),$ where $\alpha\in\mathcal {C}_1(\mathbb{A}_{2N};n,
2N)\cap\mathcal {C}_2(\mathbb{A}_{2N};n, 2N)$ and $\beta$ is a
partition with parts divisible by $2N.$ Then we insert $\beta$ into
$\alpha$ to form a partition $\lambda\in\mathcal
{C}_1(\mathbb{A}_{2N};n, 2N).$

Formally speaking, the inverse map $\Phi^{-1}$ consists of the
following two steps.

{\noindent  Step 1:} Extraction of parts from $\gamma$.

{\noindent}Suppose $\gamma\in\mathcal {C}_2(\mathbb{A}_{2N};n, 2N).$
Let $\alpha=\gamma,\,\beta=\emptyset$ and $t=l(\gamma)$. We can
obtain a pair of partitions $(\alpha, \beta)$ by  the following
procedures:
\begin{itemize}
\item If $\alpha_t$ is divisible by $N,$ then there exists an
integer $i$ such that $\alpha_t-\alpha_{t+1}=i\cdot2N+r_t$, where
$0\le r_t<2N;$ \item If $\alpha_t$ is not divisible by $N,$ then
there exists an integer $i$ such that
$\alpha_t-\alpha_{t+1}=i\cdot2N+r_t$, where $0<r_t\le2N;$ \item
Subtract $i\cdot2N$ from the parts
$\alpha_1,\alpha_2,\ldots,\alpha_t;$ Rearrange these parts to
generate a new partition $\alpha$ and add $i$ parts of size
$t\cdot2N$ to $\beta.$ \item If $t\geq 2$, then replace $t$ by
$t-1$ and repeat the above procedure. If $t=1$, we get a pair of
partitions $(\alpha,\beta)$.

\end{itemize}

{\noindent  Step 2:} Insert $\beta$ into $\alpha$.

 Assume that $(\alpha,\beta)$ is a pair of
partitions such that $\alpha\in\mathcal {C}_1(\mathbb{A}_{2N};n,
2N)\cap\mathcal {C}_2(\mathbb{A}_{2N};n, 2N)$ and $\beta$ is a
partition with parts divisible by $2N.$ We can construct a
partition $\lambda \in \mathcal {C}_1(\mathbb{A}_{2N};n, 2N)$. If
$\beta_1\le \alpha_1+2N-1,$ we insert all the parts of $\beta$
into $\alpha$ to generate a new partition $\lambda.$

 If $\beta_1 > \alpha_1+2N-1$, we set $t=1$ initially and iterate the
  following procedure until $\beta_t\leq \alpha_1+2N-1$:
\begin{itemize}
\item Let $i$ be the largest positive integer such that
$\beta_t-i\cdot2N\ge\alpha_i,$ namely for $j>i$ we have
$\beta_t-j\cdot2N<\alpha_j.$ \item Add $2N$ to the first $i$ parts
$\alpha_1,\alpha_2,\ldots,\alpha_i,$  then insert $\beta_t-i\cdot2N$
into $\alpha$ in the position before the part $\alpha_{i+1}.$ \item
Rearrange the resulted parts to form a new partition $\alpha$ and
replace $t$ by $t+1.$
\end{itemize}

Finally, we arrive at the condition $\beta_t\le\alpha_1+2N-1.$
Then we  insert all the remaining parts of $\beta$ into $\alpha$
to generate a new partition $\lambda.$ It can be shown that
 $\lambda \in \mathcal {C}_1(\mathbb{A}_{2N};n, 2N)$. For the details of the proof, see
\cite{bes91,yee02}.

 We now turn to the generalization
 of Bessenrodt's insertion algorithm and we will show how
 one can derive Theorem \ref{cgeneral} from this generalized algorithm.
Let $\mathbb{A}_N=\mathbb{A'}_N\cup\mathbb{A''}_N$ with
$\mathbb{A'}_N\cap\mathbb{A''}_N=\emptyset.$

\begin{defi}
Let $\mathcal {B}_1(\mathbb{A'}_N,\mathbb{A''}_N;n, N)$ denote the
set of partitions of $n$ satisfying the following conditions:
\begin{enumerate}
\item Each part is congruent to 0 or some $a_i$ modulo $N$;

\item Only the part congruent to 0 or  some $a_i$ belonging to
$\mathbb{A'}_N$ modulo $N$ can be repeated;

\item The difference between two successive parts is at most $N$
and strictly less than $N$ if either part congruent to 0 or some
$a_i$ belonging to $\mathbb{A'}_N$ modulo $N$;

\item The smallest part is less than $N$.
\end{enumerate}
\end{defi}

\begin{defi}
Let $\mathcal {B}_2(\mathbb{A'}_N, \mathbb{A''}_N;n, N)$ denote
the set of partitions of $n$ satisfying the following conditions:
\begin{enumerate}
\item Each part is congruent to some $a_i$ modulo $N;$

\item Only part congruent to some $a_i$ belonging to
$\mathbb{A'}_N$ modulo $N$ can be repeated.
\end{enumerate}
\end{defi}

The cardinalities of $\mathcal {B}_1(\mathbb{A'}_N,
\mathbb{A''}_N;n)$ and $\mathcal {B}_2(\mathbb{A'}_N,
\mathbb{A''}_N;n)$ are denoted by
 $b_1(\mathbb{A'
 }_N, \mathbb{A''}_N;n)$ and $b_2(\mathbb{A'}_N,\mathbb{A''}_N;n, N)$ respectively.

Bessenrodt's generalization of the Andrews-Olsson theorem is stated
as follows.

\begin{thm}[Bessenrodt]\label{bess}
For any $n \in\mathbb{N}$, we have $$ b_1(\mathbb{A'
 }_N, \mathbb{A''}_N;n,N)=b_2(\mathbb{A'}_N,\mathbb{A''}_N;n, N).$$
\end{thm}

Clearly, Andrews-Olsson's Theorem \ref{and-ols} can be viewed as
the special
case $\mathbb{A'
 }_N=\emptyset$ of Theorem \ref{bess}. Theorem \ref{cgeneral} is
 the special case for  $2N$ and $\mathbb{A'}_{2N}=\{N\}$.
 As noted by Bessenrodt \cite{bes95},  the special case $N=2$,
  $\mathbb{A'}_2=\{1\}$
  and $\mathbb{A''}_2=\emptyset$,  reduces to
  Euler's partition Theorem,  and Bessenrodt's insertion
  algorithm for
  this case  coincides with Sylvester's bijection.

\setcounter{equation}{0}
\section{Connection to Boulet's formula}

In this section, we show that our two-parameter refinement
\eqref{main2-e} can be derived from a formula of Boulet. The
following four-parameter weight was introduced by Boulet \cite{bou}
as a generalization of the weight defined by Andrews \cite{and041}.
Let $a,\,b,\,c$ and $d$ be commuting indeterminants. Define the
following weight function $\omega(\lambda)$ on the set of
partitions:
\[\omega(\lambda)=a^{\sum_{i\geq1}\lceil\lambda_{2i-1}/2\rceil}
b^{\sum_{i\geq1}\lfloor\lambda_{2i-1}/2\rfloor}
c^{\sum_{i\geq1}\lceil\lambda_{2i}/2\rceil}
d^{\sum_{i\geq1}\lfloor\lambda_{2i}/2\rfloor},\]
 where $\lceil x
\rceil$ (resp. $\lfloor x \rfloor$) stands for the smallest (resp.
largest) integer greater (resp. less) than or equal to $x$ for a
given real number $x$. Boulet obtained the following formula:
\begin{equation}\label{bous}
\sum_{\lambda \in P}\omega(\lambda)=\prod_{j=1}^{\infty}
\frac{(1+a^jb^{j-1}c^{j-1}d^{j-1})(1+a^jb^jc^jd^{j-1})}{(1-a^jb^jc^jd^j)
(1-a^jb^jc^{j-1}d^{j-1})(1-a^jb^{j-1}c^jd^{j-1})},
\end{equation}
where $P$ denotes the set of integer partitions. It can be easily
checked that the generating function of  partitions in which every
part appears an even number of times is
\[\prod_{j=1}^{\infty}\frac{1}{(1-a^jb^jc^jd^j)(1-a^jb^{j-1}c^jd^{j-1})}.\]
From \eqref{bous}, Boulet deduced the generating function for the
weight function $\omega(\lambda)$ when $\lambda$ runs over
partitions with distinct parts (\cite[Corollary 2]{bou}):
\begin{equation}\label{bous1}
\sum_{\lambda \in \mathcal{D}}\omega(\lambda)=\prod_{j=1}^{\infty}
\frac{(1+a^jb^{j-1}c^{j-1}d^{j-1})(1+a^jb^jc^jd^{j-1})}{
(1-a^jb^jc^{j-1}d^{j-1})}.
\end{equation}
Making the substitutions
 $a\mapsto xyq,\,b\mapsto x^{-1}yq,\,c\mapsto xy^{-1}q,\,
d\mapsto x^{-1}y^{-1}q$ in \eqref{bous}, Boulet derived the
following identity due to   Andrews \cite{and041}.

\begin{thm}[Andrews] We have
\[
\sum_{\lambda\in
P}x^{l_o(\lambda)}y^{l_o({\lambda}')}q^{|\lambda|}=\prod_{j=1}^{\infty}\frac{(1+xyq^{2j-1})}{(1-q^{4j})(1-x^2q^{4j-2})(1-y^2q^{4j-2})}
.\]
\end{thm}

 Using the same substitution in
\eqref{bous1}, we find obtain the following formula for partitions
with distinct parts.

\begin{thm}\label{aa} We have
\begin{equation} \label{d}
\sum_{\lambda \in
\mathcal{D}}x^{l_o(\lambda)}y^{l_o(\lambda')}q^{|\lambda|}=
\prod_{j=1}^{\infty}\frac{1+xyq^{2j-1}}{1-y^2q^{4j-2}}.
\end{equation}
\end{thm}

On the other hand, it is easy to derive the following generating
function.

\begin{thm} \label{ab}
We have
\begin{equation} \label{o}
\sum_{\lambda \in
\mathcal{O}}x^{n_o(\lambda)}y^{l(\lambda)}q^{|\lambda|}=
\prod_{j=1}^{\infty}\frac{1+xyq^{2j-1}}{1-y^2q^{4j-2}}.
\end{equation}
\end{thm}

\pf We have
\begin{align*}
\prod_{j=1}^{\infty}\frac{1+xyq^{2j-1}}{1-y^2q^{4j-2}}
&=\prod_{j=1}^{\infty}\left(\left(1+xyq^{2j-1}\right)
\sum_{i=0}^{\infty}y^{2i}q^{(2i)\cdot(2j-1)}\right)\\
&=\prod_{j=1}^{\infty}\sum_{i=0}^{\infty}\left(
y^{2i}q^{(2i)\cdot(2j-1)}+xy^{(2i+1)}q^{(2i+1)\cdot(2j-1)}\right)\\
&=\prod_{j=1}^{\infty}\left(1+xyq^{(2j-1)}+y^{2}q^{2\cdot(2j-1)}
+xy^3q^{3\cdot(2j-1)}+\cdots\right)\\
&=\sum_{\lambda \in
\mathcal{O}}x^{n_o(\lambda)}y^{l(\lambda)}q^{|\lambda|},
\end{align*}
as desired. \qed

Since $l_a(\lambda)=l_o(\lambda')$ for any
 partition $\lambda$, combining Theorems \ref{aa} and \ref{ab}
 yields Theorem \ref{main2}.

\section{A combinatorial proof of the main result}

In this section, we  give a combinatorial proof of Theorem
\ref{main2}. We will use a restricted version of the variant of
Bessenrodt's insertion algorithm given in Section 3. We now
proceed to give the proofs of Lemma \ref{lem1},  Lemma \ref{lem2}
and Theorem \ref{temp-t}.

\noindent{\it Proof of Lemma \ref{lem1}.} For $\lambda \in
\mathcal{D}(n)$, define $\alpha=\varphi(\lambda)$ to be the
$2$-modular diagram conjugate of $\lambda$. It is necessary to show
that $\alpha \in \mathcal{A}_1(n)$. On the one hand, it is easy to
see
 that there is
no
 odd part in $\alpha$ that
 can be repeated, since there is at most one
``1" in each row of the $2$-modular diagram of $\lambda$.
Moreover, the condition that $\lambda$ is a
 partition with distinct parts implies that
 the difference between successive parts in $\alpha$ is at most $4$
and strictly less than $4$ if either part is divisible by  $2$,
and that the smallest part of $\alpha$ is less than $4$.

The reverse map $\varphi^{-1}$ can be easily constructed. For
$\alpha \in \mathcal{A}_1(n)$, we note that its $2$-modular
diagram conjugate is a partition with distinct parts, namely,
$\lambda=\varphi^{-1}(\alpha)\in \mathcal{D}(n)$. Thus $\varphi$
is a bijection. Furthermore, it is not difficult to check
$l_o(\lambda)=l_o(\alpha)$ and
$l_a(\lambda)=2r_2(\alpha)+l_o(\alpha)$. This completes the proof.
\qed

\noindent{\it Proof of Lemma \ref{lem2}.} Let
$\mu=(1^{m_1},3^{m_3},5^{m_5},\ldots,(2t-1)^{m_{2t-1}}) \in
\mathcal{O}(n)$ . For every multiplicity $m_i$, we write
$m_i=2h_i+s_i$ ($s_i=0,\,1$). Then we define
$\beta=\psi(\mu)=(1^{m'_1},2^{m'_2},3^{m'_3},\ldots,k^{m'_{k}})$,
where $m'_{2i+1}=s_{2i+1}$ and $m'_{2i}=h_{i}$. Clearly,
$m'_{4i}=h_{2i}=0$ and $m'_{2i+1}\leq 1$, and so $\beta \in
\mathcal{A}_2(n).$ For example, let $\mu=(1,3,7^2,9,15)$. Then we
have $\beta=(1,3,9,14,15)$ whose  2-modular diagram is illustrated
in Figure \ref{fig1}.

\begin{figure}[h]
\begin{center}
\setlength{\unitlength}{0.5mm}
\begin{picture}(90,50)
\put(10,50){\line(0,-1){50}} \put(20,50){\line(0,-1){50}}
\put(30,50){\line(0,-1){40}} \put(40,50){\line(0,-1){30}}
\put(50,50){\line(0,-1){30}} \put(60,50){\line(0,-1){30}}
\put(70,50){\line(0,-1){20}} \put(80,50){\line(0,-1){20}}
\put(90,50){\line(0,-1){10}} \put(10,0){\line(1,0){10}}
\put(10,10){\line(1,0){20}} \put(10,20){\line(1,0){50}}
\put(10,30){\line(1,0){70}} \put(10,40){\line(1,0){80}}
\put(10,50){\line(1,0){80}} \multiput(13,42)(10,0){7}{$2$}
\multiput(13,32)(10,0){7}{$2$} \multiput(13,22)(10,0){4}{$2$}
\put(83,42){$1$} \put(53,22){$1$} \put(13,12){$2$}
\put(23,12){$1$} \put(13,2){$1$} \put(0,43){$15$}
\put(0,33){$7^2$} \put(0,23){$9$} \put(0,13){$3$} \put(0,3){$1$}
\end{picture}

\caption{The diagram representation of
$\beta=(1,3,7^2,9,15)$.}\label{fig1}
\end{center}

\end{figure}
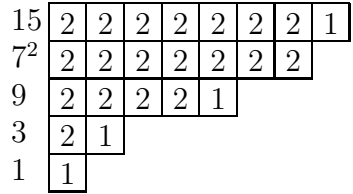

The inverse map $\psi^{-1}$ can be  easily described.  Let
$\beta=(1^{m_1},2^{m_2},3^{m_3},\ldots,t^{m_t}) \in
\mathcal{A}_2(n)$. Then we have  $m_{4i}=0$ and $m_{2i-1}=0$ or $1$
for $i\geq 1$. Let
\[  \mu =\psi^{-1}(\beta)=(1^{m'_1},2^{m'_2},3^{m'_3},\ldots,
k^{m'_k}),\] where $m'_{2i-1}=2m_{4i-2}+m_{2i-1}$ and $m'_{2i}=0$
for $i=1,2,\ldots$. Obviously, $\psi^{-1}(\beta) \in
\mathcal{O}(n)$. It follows that
 $n_o(\mu)=l_o(\beta)$ and $l(\mu)=2r_2(\beta)+l_o(\beta)$. This completes
 the proof.
\qed

It is easy to see that  Theorem \ref{temp-t} follows from
 Theorem \ref{cgeneral} by  setting $N=2$ and
 $\mathbb{A}_4=\{1,2,3\}$, that is,
$\mathcal{C}_1(\{1,2,3\};n,4)=\mathcal{A}_1(n)$ and
$\mathcal{C}_2(\{1,2,3\};n,4)=\mathcal{A}_2(n)$.

Figures \ref{fig2} and \ref{fig3} illustrate the procedure in
Theorem \ref{temp-t}.

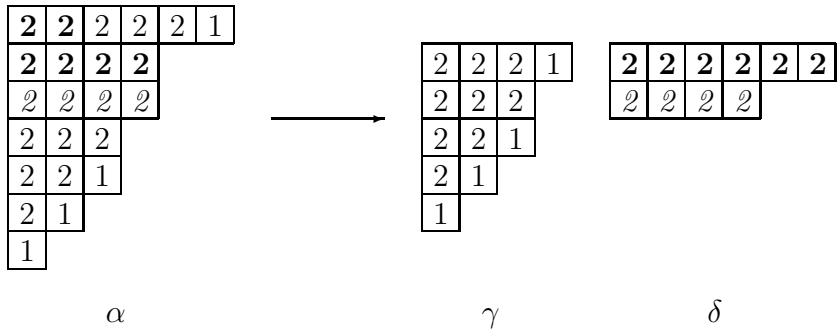
\begin{figure}[h]
\begin{center}
\setlength{\unitlength}{1mm} \psset{unit=1.2mm}
\begin{picture}(120,45)

\put(13,3){$\alpha$}
\put(0,10){\line(0,1){35}}\put(5,10){\line(0,1){35}}
\put(10,15){\line(0,1){30}}\put(15,20){\line(0,1){25}}
\put(20,30){\line(0,1){15}}\put(25,40){\line(0,1){5}}
\put(30,40){\line(0,1){5}} \put(0,10){\line(1,0){5}}
\put(0,15){\line(1,0){10}}\put(0,20){\line(1,0){15}}
\put(0,25){\line(1,0){15}} \put(0,30){\line(1,0){20}}
\put(0,35){\line(1,0){20}}\put(0,40){\line(1,0){30}}
\put(0,45){\line(1,0){30}} \multiput(1.5,41)(5,0){2}{$\textbf{2}$}
\multiput(11.5,41)(5,0){3}{$2$}
\multiput(1.5,36)(5,0){4}{$\textbf{2}$}
\multiput(1.5,26)(5,0){3}{$2$}\multiput(1.5,21)(5,0){2}{$2$}
\put(26.5,41){$1$}\put(11.5,21){$1$}
\put(1.5,16){$2$}\put(6.5,16){$1$}
\put(1.5,11){$1$}\multiput(1.5,31)(5,0){4}{$\emph{2}$}

\put(63,3){$\gamma$} \put(35,30){\vector(1,0){15}}
\put(55,15){\line(0,1){25}}\put(60,15){\line(0,1){25}}
\put(65,20){\line(0,1){20}}\put(70,25){\line(0,1){15}}
\put(75,35){\line(0,1){5}}
\put(55,15){\line(1,0){5}}\put(55,20){\line(1,0){10}}
\put(55,25){\line(1,0){15}}\put(55,30){\line(1,0){15}}
\put(55,35){\line(1,0){20}} \put(55,40){\line(1,0){20}}
\multiput(56.5,36)(5,0){3}{$2$}\multiput(56.5,31)(5,0){3}{$2$}
\multiput(56.5,26)(5,0){2}{$2$}\multiput(56.5,16)(5,5){3}{$1$}
\put(71.5,36){$1$}\put(56.5,21){$2$}

\put(93,3){$\delta$}
\put(80,30){\line(1,0){20}}\put(80,35){\line(1,0){30}}
\put(80,40){\line(1,0){30}}
\multiput(80,30)(5,0){5}{\line(0,1){10}}
\multiput(105,35)(5,0){2}{\line(0,1){5}}
\multiput(81.5,36)(5,0){6}{$\textbf{2}$}\multiput(81.5,31)(5,0){4}{$\emph{2}$}
\end{picture}

\caption{Extraction of parts from $\alpha$.}\label{fig2}
\end{center}
\end{figure}

\begin{figure}[h]
\begin{center}\setlength{\unitlength}{1mm}
\begin{picture}(120,45)
\put(8,5){$\gamma$}
\put(0,15){\line(0,1){25}}\put(5,15){\line(0,1){25}}
\put(10,20){\line(0,1){20}}\put(15,25){\line(0,1){15}}
\put(20,35){\line(0,1){5}}
\put(0,15){\line(1,0){5}}\put(0,20){\line(1,0){10}}
\put(0,25){\line(1,0){15}}\put(0,30){\line(1,0){15}}
\put(0,35){\line(1,0){20}} \put(0,40){\line(1,0){20}}
\multiput(1.5,36)(5,0){3}{$2$}\multiput(1.5,31)(5,0){3}{$2$}
\multiput(1.5,26)(5,0){2}{$2$}\multiput(1.5,16)(5,5){3}{$1$}
\put(16.5,36){$1$}\put(1.5,21){$2$}

\put(33,5){$\delta$}
\put(25,30){\line(1,0){20}}\put(25,35){\line(1,0){30}}
\put(25,40){\line(1,0){30}}
\multiput(25,30)(5,0){5}{\line(0,1){10}}
\multiput(50,35)(5,0){2}{\line(0,1){5}}
\multiput(26.5,36)(5,0){6}{$\textbf{2}$}\multiput(26.5,31)(5,0){4}{$\emph{2}$}
\put(60,32){\vector(1,0){15}}
\put(80,15){\line(0,1){25}}\put(85,15){\line(0,1){25}}
\put(90,20){\line(0,1){20}}\put(95,25){\line(0,1){15}}
\put(100,25){\line(0,1){15}}\put(105,25){\line(0,1){15}}
\put(110,30){\line(0,1){10}}\put(115,30){\line(0,1){10}}
\put(120,35){\line(0,1){5}}

\put(95,5){$\beta$}
\put(80,15){\line(1,0){5}}\put(80,20){\line(1,0){10}}
\put(80,25){\line(1,0){25}}\put(80,30){\line(1,0){35}}
\put(80,35){\line(1,0){40}}\put(80,40){\line(1,0){40}}

\multiput(81.5,21)(0,5){4}{$2$}\multiput(86.5,26)(0,5){3}{$2$}
\multiput(91.5,26)(0,5){3}{$\textbf{2}$}\multiput(96.5,26)(0,5){3}{$\textbf{2}$}
\multiput(101.5,31)(0,5){2}{$2$}\multiput(106.5,31)(0,5){2}{$\emph{2}$}
\multiput(111.5,31)(0,5){2}{$\emph{2}$}
\put(116.5,36){1}\put(116.5,36){1} \put(101.5,26){1}
\put(86.5,21){1}\put(81.5,16){1}

\end{picture}

\caption{Insertion all parts of $\delta$ into
$\gamma$.}\label{fig3}
\end{center}
\end{figure}
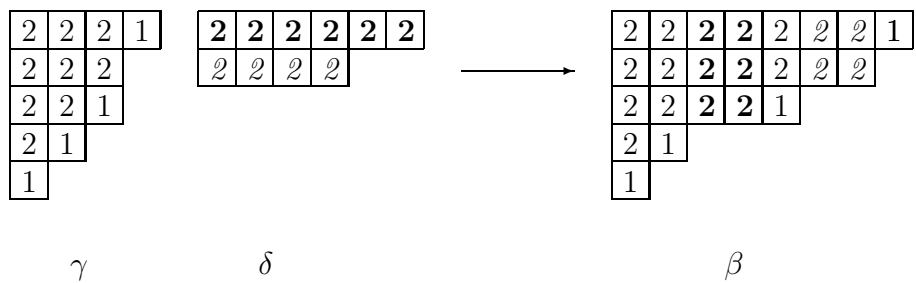

To conclude, we combine the maps $\psi$, $\phi$ and $\varphi$ to
construct the desired map $\Delta$  from the set of partitions with
distinct parts to the set of partitions with odd parts parts:
 $\Delta=\psi^{-1}\circ\phi\circ\varphi$. The properties of
 $\Delta$ lead to a proof of Theorem \ref{main2}.

For example, let $\lambda=(17,16,14,10,7,4,2,1) \in
\mathcal{D}(71).$ Then we have
\begin{eqnarray*}
 \alpha & = & \varphi(\lambda)=(15,12,10, 9,8,6,6,4,1),\\[3pt]
\beta & = & \phi(\alpha)=(19,18,13,10,6,5),\\[3pt]
 \mu & =  & \psi^{-1}(\beta)=(3^2,5^3,9^2,13,19) \in \mathcal{O}(71).
 \end{eqnarray*}
 Moreover,
$l_o(\lambda)=3$, $l_a(\lambda)=9$, $n_o(\mu)=3$ and $l(\mu)=9$.

\vspace{.2cm}
\noindent{\bf Acknowledgments.}  This
work was supported by the 973 Project, the PCSIRT Project of the
Ministry of Education, the Ministry of Science and Technology,  and
the National Science Foundation of China.

\end{document}